\documentclass{elsart}
\usepackage{ifpdf}
\usepackage{graphicx,amssymb,lineno}
\ifpdf
\usepackage[%
  pdftitle={Instructions for use of the document class
    elsart},%
  pdfauthor={Simon Pepping},%
  pdfsubject={The preprint document class elsart},%
  pdfkeywords={instructions for use, elsart, document class},%
  pdfstartview=FitH,%
  bookmarks=true,%
  bookmarksopen=true,%
  breaklinks=true,%
  colorlinks=true,%
  linkcolor=blue,anchorcolor=blue,%
  citecolor=blue,filecolor=blue,%
  menucolor=blue,pagecolor=blue,%
  urlcolor=blue]{hyperref}
\else
\usepackage[%
  breaklinks=true,%
  colorlinks=true,%
  linkcolor=blue,anchorcolor=blue,%
  citecolor=blue,filecolor=blue,%
  menucolor=blue,pagecolor=blue,%
  urlcolor=blue]{hyperref}
\fi

\makeatletter
\def\elsartstyle{%
    \def\normalsize{\@setfontsize\normalsize\@xiipt{14.5}}
    \def\small{\@setfontsize\small\@xipt{13.6}}
    \let\footnotesize=\small
    \def\large{\@setfontsize\large\@xivpt{18}}
    \def\Large{\@setfontsize\Large\@xviipt{22}}
    \skip\@mpfootins = 18\p@ \@plus 2\p@
    \normalsize
}
\@ifundefined{square}{}{}
\makeatother

\def\be{\begin{equation}}
\def\ee{\end{equation}}

\begin{document}
\begin{frontmatter}
\title{Two-dimensional Systems that Arise from the Noether Classification of
Lagrangians on the Line}
\author{M. Umar Farooq$^{a}$, S. Ali$^{b}$ and F. M. Mahomed$^{c}$}

\address{$^{a}${Center for Advanced Mathematics and Physics, 
National University of Sciences and Technology, H-12 Campus, Islamabad 44000, Pakistan.}\\
$^{b}${School of Electrical Engineering and Computer Science, 
National University of Sciences and Technology,  H-12 Campus, Islamabad 44000, Pakistan.}\\
$^{c}${Centre for Differential Equations, Continuum Mechanics and
Applications, School of Computational and Applied Mathematics,
University of the Witwatersrand, Wits 2050, South Africa.}}

\begin{abstract}
Noether-like operators play an essential role in writing down the
first integrals for Euler-Lagrange systems of ordinary differential
equations (ODEs). The classification of such operators is carried
out with the help of analytic continuation of Lagrangians on the
line. We obtain the classification of 5, 6 and 9 Noether-like
operators for two-dimensional Lagrangian systems that arise from the
submaximal and maximal dimensional Noether point symmetry
classification of Lagrangians on the line. Cases in which the
Noether-like operators are also Noether point symmetries for the
systems of two ODEs are  mentioned. In particular, the
8$-$dimensional maximal Noether algebra is remarkably obtained for
the simplest system of the free particle equations in two dimensions
from the 5$-$dimensional complex Noether algebra of the standard
Lagrangian of the scalar free particle equation. We present the
effectiveness of Noether-like operators for the determination of
first integrals of systems of two nonlinear differential equations
which arise from scalar complex Euler-Lagrange ODEs that admit
Noether symmetry. 
\end{abstract}


\end{frontmatter}

\section{Introduction}
\qquad Over the last few decades there have been many contributions to the study of
Lie and Noether point symmetries for second-order ODEs. The use of these
symmetries signifies their importance in the reduction of order of a given
dynamical equation and in constructing its first integrals (constants of the
motion). Lie proved that the maximum dimension of the point symmetry algebra
for a scalar second-order ODE is eight. He showed that the simple differential equation $%
u^{\prime \prime }=0,$ admits the maximal symmetry algebra which corresponds
to the Lie algebra $sl(3,\Re)$ and there exists a linearizing point transformation for those second-order
scalar ODEs that possess this algebra. Further, any other algebra admitted
by a scalar second-order ODE is a subalgebra of $sl(3,\Re)$ (see  \cite{1}).
Lie also proved that a scalar second-order ODE admits the
maximal $\mathit{r}\in \{0,1,2,3,8\}-$dimensional point symmetry algebras
\cite{2} (see also \cite{1}). Moreover, it is also known that the maximum
Lie algebra of point symmetries admitted by a system of two second-order
ODEs is fifteen \cite{3}.

There are at least two important facets in Lagrangian mechanics. One
is to find a Lagrangian and the other is to construct first integrals of
the underlying Euler-Lagrange (EL) differential equations.
The classical Noether theorem requires the presence
of a Lagrangian for a differential equation before its first integral can
be evaluated with the help of an explicit formula \cite{4}. The
 first integrals not only play a vital role in the integrability of an
ODE but also have direct physical consequences. Douglas \cite{5} gave the
solution to the inverse problem for a system of two second-order ODEs. It is
at times very difficult to find a Lagrangian and there are differential
equations (see, e.g., Anderson and Thompson
\cite{6}) that do not admit Lagrangians. It raises another important
question: can we find first integrals without a variational structure? The
construction of conservation laws in the absence of a Lagrangian is carried
out by Kara and Mahomed \cite{7}. Ibragimov \cite{8}, also discussed a
way of finding conservation laws without prior knowledge of a Lagrangian.
In \cite{9}, Kara and Mahomed introduced the partial
Noether approach to construct conserved quantities. They derived first
integrals for those differential equations that fail to admit Lagrangians. Furthermore 
for nonvariational differential equations the Noether appraoch is also used 
to find first integrals recently by Gouveia and Torres in \cite{10}. Moreover 
the proof of more general forms of the Noether theorem and conservation
laws is a subject under strong development see, e.g., \cite{11,12,13}. 

The search for new conservation laws for those systems of two ODEs that
appear as dynamical equations in mathematical physics is indispensable. In
\cite{14}, Gorringe and Leach classified the Lie algebra for systems of two
second-order ODEs with constant coefficients. Later, Wafo and Mahomed
discussed the case with variable coefficients in \cite{15}. The
classification of Noether point symmetries of an autonomous quadratic in the velocities Lagrangian
with two degrees of
freedom is done by Sen \cite{16}. We know that the maximum dimension of the
Noether point symmetry algebra for a Lagrangian in one-dimensional particle
dynamics is five. A first-order Lagrangian on the line can have Noether
point symmetry algebras of dimension $0,1,2,3$ or $5$ \cite{17,18,19}.
Besides, a two dimensional system of free particle equations has $8-$%
dimensional Noether algebra \cite{3,20}. The classification of
Noether subalgebras for systems of two nonlinear EL ODEs
has not been carried out before. Our aim is to classify Noether
operators for two-dimensional systems that arise from the complex
Noether point symmetry classification of complex Lagrangians on the
line. We make use of Noether-like operators for such nonlinear
systems by utilizing complex arguments.

In \cite{21,22,23}, the authors used analytic continuation of ODEs
in the complex plane to obtain non-trivial results for systems
of differential equations. In particular, symmetry analysis on the complex plane leads
to finding the symmetries, reduction of order, linearization and
conservation laws of systems of two real ODEs and PDEs. In this
paper, we furnish important results on Noether operators and the
production of first integrals of systems of two second-order ODEs.
It is mainly done by confining the domain of definition of dependent
functions on a single real line and then the dynamics is studied on
the complex plane. We employ a procedure which is similar to
analytic continuation in the intervening steps but require that the
dependent variables may be restricted on the line to achieve
operators and first integrals of systems of ODEs. To achieve our
goal we start with restricted complex ordinary differential
equations (r-CODEs) \cite{21}. Such r-CODEs are obtained by allowing
a complex function to depend only on a single real variable which
thereby yields a system of two\ ODEs. Our aim is thus to investigate
the algebraic  properties of systems of EL ODEs
encoded in r-CODEs that admit Lagrangians. Indeed, it is shown that
a complex Lagrangian reveals two real inequivalent Lagrangians for
systems of two ODEs, i.e., these do not differ by a divergence.
Thus, the classification of Lagrangians on the line would offer us a
great deal of information about the inverse problem, algebraic
properties and first integrals of the corresponding systems of two
ODEs. In this regard, we mention that the complex extension of
Hamiltonians and Lagrangians has been discussed by Bender \cite{24}.
He proved various intriguing results that verified known facts about
quantum mechanics. The essence of his approach lies in the fact that
a complex Hamiltonian can also be taken for a consistent physical
theory of quantum mechanics.

Unlike a complex Lagrangian, a complex symmetry of an r-CODE may not give
two real symmetries of the system, in general. It splits into two operators.
We call such operators  \textit{Lie-like operators}. These operators may
not necessarily be symmetries of the systems of ODEs corresponding to an
r-CODE. Further, these operators do not form an algebra in general.
Similarly, we call \textit{Noether-like operators} those operators that are obtainable from a
complex Noether symmetry. To draw attention on the significance of these
operators, a few examples are discussed in some detail. The complex Noether theorem enables
a description of an explicit formula for Noether-like operators. It is
shown explicitly that the first integrals for such systems of two
second-order ODEs corresponding to these operators associated with the
Lagrangians can be determined. An appealing consequence emerges when a single
real symmetry permits the existence of two first integrals for the corresponding system.

We obtain a classification of a system of two EL ODEs with respect
to the Noether-like operators they admit. We see that systems of
ODEs that exhibit the same invariance properties as r-CODEs have the
same structure of the Noether-like operators. We construct an
analogue of the Noether counting theorem for systems of two
second-order EL ODEs that arise from the Lagrangian formulation of
scalar second-order equations. We obtain $5$ and $6$ Noether-like
operators for systems of two second-order ODEs by the analytic
continuation of the $3-$dimensional Noether algebra of an r-CODE.
The Noether algebra of a restricted complexified free particle
equation which is 5-dimensional remarkably generates an
$8-$dimensional Noether algebra of the simplest system which
exhibits a non-trivial implication of the complex variable approach.
Moreover, \textit{nine} Noether-like operators imply \textit{ten}
first integrals whereas in the classical Noether approach, there are
only eight first integrals corresponding to eight Noether
symmetries. This indeed is a nice result which we derive here. It is
also conjectured that all linear systems of two second-order ODEs
that can be derived from a complex variational principal of a scalar
linear second-order ODE, indeed, admit an eight dimensional Noether
algebra.

The outline of the paper is as follows. In the next section, we
present the preliminaries. The EL equations are obtained and the
Noether theorem is invoked in order to write down the Noether-like
symmetry conditions for systems of two ODEs. An explicit formula to
find first integrals corresponding to Noether-like operators is also
mentioned for such systems. Few cases are described that shed light
on the comparison of two formulae for acquiring first integrals for
systems of ODEs. The classification of $5$ and $6$ Noether-like
operators is carried out in the third section. The second last
section deals with the case of maximal Noether symmetry algebra that
is attained via complex variables. Some physical insights are also
developed in the same section. Finally, we conclude the discussion
in the last section.

\section{Preliminaries}
The problem of our interest is to study and classify algebraic properties and invariants for systems
of two second-order ODEs of the form%
\begin{eqnarray}
f^{\prime \prime } =w_{1}(x,f,g,f^{\prime },g^{\prime }),  \nonumber \\
g^{\prime \prime } =w_{2}(x,f,g,f^{\prime },g^{\prime }), \label{eq1}
\end{eqnarray}%
by invoking a complex analytic structure on the $fg$-plane. This
study helps in discovering conserved quantities for those systems of
ODEs that appear in diverse physical phenomena. For example, in the
case of two coupled nonlinear oscillators the governing dynamical
equations have the form of (\ref{eq1}). There has been  considerable amount
of work done in finding all the conserved quantities of time
dependent and independent nonlinear oscillators. The search for new
conservation laws for nonlinear systems is one of the main
objectives of the physicists. We assume a variational structure on
the system (\ref{eq1}). We show how the Noether point symmetry
classification of complex Lagrangians on the line guarantees the
emergence of operators and invariants for a class of systems of the
form (\ref{eq1}). It is emphasized that the use of complex Lagrangians is
inevitable for such systems.

We assume that the above system has a complex structure, i.e., there
exists a transformation
\begin{equation}
u(x)=f(x)+ig(x),\label{eq2}
\end{equation}%
that maps system (\ref{eq1}) to
\begin{equation}
u^{\prime \prime }=w(x,u,u^{\prime }), \label{eq3}
\end{equation}%
which is a second-order r-CODE. It is clear that an arbitrary system
(\ref{eq1}) may or may not necessarily correspond to an r-CODE. Here we
are interested in looking for the insights that can be extracted for
systems of two ODEs (\ref{eq1}) from equations of the form (\ref{eq3}). The
symmetry analysis of such systems $(1)$ is carried out in
\cite{21,22,23}. The complex Lie algebra of (\ref{eq3}) gives the real
Lie-like operators of the system $(1)$ \cite{21,22,23}. We intend to
investigate the implications of the variational structure on such
systems of differential equations via the complex plane. The
following details are not in conflict with conventional symmetry
analysis but is rather a natural generalization of it in the complex
domain. The equation (\ref{eq3}) may be regarded as an analytic
continuation of a general scalar second-order ODE in the restricted
complex domain, i.e., here $u$ is a complex function of a real
variable $x$.

Suppose that the equation (\ref{eq3}) appears from a variational
principle, i.e., there exists a complex Lagrangian $L(x,u,u^{\prime
})$ such that the EL equation implies (\ref{eq3}). Once the Lagrangian of
a differential equation is known, its symmetry properties are
explored. In the subsequent discussions we consider EL equations,
Noether-like symmetry
conditions and formulae of first integrals for the systems of two ODEs.%
\newline
\textbf{Theorem. }\textit{If }$L(x,u,u^{\prime })=L_{1}+iL_{2},$\textit{\ is
a complex Lagrangian of a second-order r-CODE then both }$L_{1}(x,f,g,f^{\prime
},g^{\prime })$\textit{\ and }$L_{2}(x,f,g,f^{\prime },g^{\prime })$\textit{%
\ are two Lagrangians of the corresponding system of two second-order ODEs.}\newline
\textbf{Proof. }\ Suppose that $L$ is a complex Lagrangian of the r-CODE (\ref{eq3}) relative
to system (\ref{eq1}). Therefore, it satisfies the complex EL equation. The realification of the EL
equation yields%
\begin{eqnarray}
\frac{\partial L_{1}}{\partial f}+\frac{\partial L_{2}}{\partial g}-\frac{d}{%
dx}\left(\frac{\partial L_{1}}{\partial f^{\prime }}+\frac{\partial L_{2}}{%
\partial g^{\prime }}\right) =0,  \nonumber \\
\frac{\partial L_{2}}{\partial f}-\frac{\partial L_{1}}{\partial g}-\frac{d}{%
dx}\left(\frac{\partial L_{2}}{\partial f^{\prime }}-\frac{\partial L_{1}}{%
\partial g^{\prime }}\right) =0. \label{eq4}
\end{eqnarray}%
Since $L(x,u,u^{\prime })$ is complex analytic in its arguments,
both $L_{1}$ and $L_{2}$ satisfy the Cauchy-Riemann equations and the above
system becomes%
\begin{eqnarray}
\frac{\partial L_{1}}{\partial f}-\frac{d}{dx}\left(\frac{\partial L_{1}}{%
\partial f^{\prime }}\right) =0,\frac{\partial L_{1}}{\partial g}-\frac{d}{dx}\left (%
\frac{\partial L_{1}}{\partial g^{\prime }}\right )=0,  \nonumber \\
\frac{\partial L_{2}}{\partial f}-\frac{d}{dx}\left (\frac{\partial L_{2}}{%
\partial f^{\prime }}\right) =0,\frac{\partial L_{2}}{\partial g}-\frac{d}{dx}\left(%
\frac{\partial L_{2}}{\partial g^{\prime }}\right)=0. \label{eq5}
\end{eqnarray}%
The above equations are the usual EL-equations for the system (\ref{eq1}). Hence $L_1$ and $L_2$ are two Lagrangians for the
system (\ref{eq1}).

Notice that these two Lagrangians do not differ by a divergence. It would be
nice to deduce the symmetry properties of these Lagrangians via the complex Lagrangian. We would
come to this point later as well as in the discussion of alternative Lagrangians.%
\newline
\textbf{Definition.} The operators $\mathbf{X}_{1}=2\varsigma _{1}\partial _{x}+\chi _{1}\partial
_{f}+\chi _{2}\partial _{g}$ and $\mathbf{X}_{2}=2\varsigma _{2}\partial _{x}+\chi _{2}\partial
_{f}-\chi _{1}\partial _{g}$ are said to
be Noether-like operators of system (\ref{eq1}) with respect to the Lagrangians $L_{1}$
and $L_{2}$ if they satisfy
\begin{eqnarray}
\mathbf{X}_{1}^{(1)}L_{1}-\mathbf{X}_{2}^{(1)}L_{2}+(d_{x}\varsigma
_{1})L_{1}-(d_{x}\varsigma _{2})L_{2}=d_{x}A_{1},\; d_x=d/dx, \nonumber \\
\mathbf{X}_{1}^{(1)}L_{2}+\mathbf{X}_{2}^{(1)}L_{1}+(d_{x}\varsigma
_{1})L_{2}+(d_{x}\varsigma _{2})L_{1}=d_{x}A_{2}, \label{eq6}
\end{eqnarray}%
for suitable functions $A_1$ and $A_2$.\newline
We now set
\begin{equation}
\varsigma  =\varsigma _{1}+i\varsigma _{2},\quad A=A_{1}+iA_{2},\quad 
\mathbf{Z} =\mathbf{X}_{1}+i\mathbf{X}_{2}.\label{eq7}
\end{equation}%
Further, if we let $\chi =\chi _{1}+i\chi _{2}$ and $\chi ^{(1)}=\chi
_{1}^{(1)}+i\chi _{2}^{(1)}$ in the complex symmetry%
\begin{equation}
\mathbf{Z}^{(1)}=\varsigma \frac{\partial }{\partial x}+\chi \frac{\partial
}{\partial u}+\chi ^{(1)}\frac{\partial }{\partial u^{\prime }}, \label{eq8}
\end{equation}
then $\mathbf{X}_{1}^{(1)}$ and $\mathbf{X}_{2}^{(1)}$ are%
\begin{eqnarray}
\mathbf{X}_{1}^{(1)} =2\varsigma _{1}\partial _{x}+\chi _{1}\partial
_{f}+\chi _{2}\partial _{g}+\chi _{1}^{(1)}\partial _{f^{\prime }}+\chi
_{2}^{(1)}\partial _{g^{\prime }},  \nonumber \\
\mathbf{X}_{2}^{(1)} =2\varsigma _{2}\partial _{x}+\chi _{2}\partial
_{f}-\chi _{1}\partial _{g}+\chi _{2}^{(1)}\partial _{f^{\prime }}-\chi
_{1}^{(1)}\partial _{g^{\prime }}. \label{eq9}
\end{eqnarray}%
We state these as the first prolongations of the Noether-like operators $\mathbf{X}_{1}$ and $\mathbf{X}_{2}$.

The above conditions are different from the usual Noether conditions for the
systems of two ODEs. These operators do not form an algebra in general. It may be
questioned as to what is the use of these operators. The point is that we
can determine invariants by employing such operators. To provide
concrete basis to our argument, we mention two cases in
which the conditions (\ref{eq6}) reduce to the usual Noether conditions: (a) if $%
\mathbf{Z}$\textbf{\ }has either a pure real or pure imaginary form then $%
\mathbf{Z}$ becomes a Noether symmetry for both Lagrangians. For example $%
\mathbf{Z}=\mathbf{X}_{1}$, i.e., it has only a real part so (\ref{eq6})
takes the form
\begin{eqnarray}
\mathbf{X}_{1}^{(1)}L_{1}+(d_{x}\varsigma _{1})L_{1} =d_{x}A_{1},  \nonumber \\
\mathbf{X}_{1}^{(1)}L_{2}+(d_{x}\varsigma _{1})L_{2} =d_{x}A_{2}, \label{eq10}
\end{eqnarray}%
and hence $\mathbf{X}_{1}$ \textit{is} a Noether symmetry surprisingly for the
system (\ref{eq1}) relative to both inequivalent Lagrangians $L_{1}$ and $L_{2}$,
(b) if $L=L_{1}$ then (\ref{eq6}) becomes
\begin{eqnarray}
\mathbf{X}_{1}^{(1)}L_{1}+(d_{x}\varsigma _{1})L_{1} =d_{x}A_{1},  \nonumber \\
\mathbf{X}_{2}^{(1)}L_{1}+(d_{x}\varsigma _{2})L_{1} =d_{x}A_{2}, \label{eq11}
\end{eqnarray}%
in which case $\mathbf{X}_{1}$ and $\mathbf{X}_{2}$ \textit{turn out }to be two
distinct Noether symmetries for the system (\ref{eq1}) corresponding to $L_{1}$.%
\newline
\textbf{Noether-like Theorem.} \textit{If }$\mathbf{X}_{1}$\textit{\ and }$%
\mathbf{X}_{2}$\textit{\ are two Noether-like operators of (\ref{eq1}) with respect to
the Lagrangians }$L_{1}$\textit{\ and }$L_{2}$\textit{\ then (\ref{eq1}) admits two
first integrals }%
\begin{eqnarray}
{\small I}_{1}={\small \varsigma }_{1}{\small L}_{1}{\small -\varsigma }_{2}%
{\small L}_{2}{\small +\partial _{f^{\prime }}L_{1}(\chi }_{1}{\small -f}%
^{\prime }{\small \varsigma }_{1}{\small -g}^{\prime }{\small \varsigma }_{2}%
{\small )-\partial _{f^{\prime }}L_{2}(\chi }_{2}{\small -f}^{\prime }%
{\small \varsigma }_{2}{\small -g}^{\prime }{\small \varsigma }_{1}{\small %
)-A}_{1}{\small ,}  \nonumber \\
{\small I}_{2}={\small \varsigma }_{1}{\small L}_{2}{\small +\varsigma }_{2}%
{\small L}_{1}{\small +\partial _{f^{\prime }}L_{2}(\chi }_{1}{\small -f}%
^{\prime }{\small \varsigma }_{1}{\small -g}^{\prime }{\small \varsigma }_{2}%
{\small )+\partial _{f^{\prime }}L_{1}(\chi }_{2}{\small -f}^{\prime }%
{\small \varsigma }_{2}{\small -g}^{\prime }{\small \varsigma }_{1}{\small %
)-A}_{2}{\small .} \label{eq12}
\end{eqnarray}%
\textbf{Proof.} That $I_{1}$ and $I_{2}$ are two first integrals of system (\ref{eq1})
 can be verified by
\begin{equation}
d_{x}I_{1}=0,~ d_{x}I_{2}=0.\label{eq13}
\end{equation}
on the system (\ref{eq1}).

Notice that formulae (\ref{eq12}) for the first integrals are different from the
usual Noether first integrals of a system. We can find the usual formulae for the
above mentioned cases (a) and (b) as follows. In case (a), we get
\begin{eqnarray}
{\small I}_{1} ={\small \varsigma }_{1}{\small L}_{1}{\small -f}^{\prime }%
{\small \varsigma }_{1}{\small (\partial }_{f^{\prime }}{\small L}_{1}%
{\small )+g}^{\prime }{\small \varsigma }_{1}{\small (\partial }_{f^{\prime
}}{\small L}_{2}{\small )-A}_{1},  \nonumber \\
{\small I}_{2} ={\small \varsigma }_{1}{\small L}_{2}{\small -f}^{\prime }%
{\small \varsigma }_{1}{\small (\partial }_{f^{\prime }}{\small L}_{2}%
{\small )-g}^{\prime }{\small \varsigma }_{1}{\small (\partial
}_{f^{\prime }}{\small L}_{1}{\small )-A}_{2}{\small .}\label{eq14}
\end{eqnarray}%
It may seem strange that we are \textit{obtaining} two first
integrals for (\ref{eq1}) corresponding to a single symmetry
$\mathbf{X}_{1}$\textbf{\ }(or $\mathbf{X}_{2} $). The point is that
$\mathbf{X}_{1}$ (or $\mathbf{X}_{2}$) is the Noether symmetry of
both $L_{1}$ and $L_{2}$. A simple question arises as to which of
the Lagrangians results in maximum number of Noether symmetries. The
answer helps in writing all the first integrals of systems of two
ODEs corresponding to such Noether symmetries. It is quite natural
to see that in case (b) we get two real first integrals
corresponding to $\mathbf{X}_{1}$ and $\mathbf{X}_{2}$ associated
with $L_{1},$ i.e.,
\begin{eqnarray}
{\small I}_{1} ={\small \varsigma }_{1}{\small L}_{1}+{\small (\chi }_{1}%
{\small -f}^{\prime }{\small \varsigma }_{1}{\small -g}^{\prime }{\small %
\varsigma }_{2}{\small )\partial }_{f^{\prime }}{\small L}_{1}+{\small (\chi
}_{2}{\small -f}^{\prime }{\small \varsigma }_{2}{\small -g}^{\prime }%
{\small \varsigma }_{1}{\small )\partial }_{g^{\prime }}{\small L}_{1}%
{\small -A}_{1}{\small ,}  \nonumber \\
{\small I}_{2} ={\small \varsigma }_{2}{\small L}_{1}-{\small (\chi }_{1}%
{\small -f}^{\prime }{\small \varsigma }_{1}{\small -g}^{\prime }{\small %
\varsigma }_{2}{\small )\partial }_{g^{\prime }}{\small L}_{1}+{\small (\chi
}_{2}{\small -f}^{\prime }{\small \varsigma }_{2}{\small -g}^{\prime }%
{\small \varsigma }_{1}{\small )\partial }_{f^{\prime }}{\small L}_{1}%
{\small -A}_{2}{\small .} \label{eq15}
\end{eqnarray}%
The two first integrals satisfy the coupled equations%
\begin{eqnarray}
\mathbf{X}_{1}^{(1)}I_{1}-\mathbf{X}_{2}^{(1)}I_{2} =0,  \nonumber \\
\mathbf{X}_{1}^{(1)}I_{2}+\mathbf{X}_{2}^{(1)}I_{1} =0, \label{eq16}
\end{eqnarray}%
where $\mathbf{Z}^{(1)}=\mathbf{X}_{1}^{(1)}+i\mathbf{X}_{2}^{(1)}.$

To illustrate we commence with those systems of second-order ODEs that admit $%
1,2,3$ or $4$ Noether-like operators. In order to carry out 
calcluations we use Computer Algebra System (CAS), e.g., MAPLE and 
CRACK \cite{25,26}. The first two examples briefly explains the determintation of 
these operators and first integrals. It is noticed that these
operators could be supplied by the analytical continuation of
$2-$dimensional Noether algebras. There exists two types of
realizations of both one and two dimensional complex algebras
\cite{2,19}. The one dimensional realization includes $\partial
/\partial x$ and $\partial /\partial u$ that correspond to a single
and two Noether-like operators. Similarly, $\{\partial /\partial
x,\partial /\partial u\}$ and $\{\partial /\partial u,x\partial
/\partial x+u\partial /\partial u\}$ are amongst the  realizations
of two dimensional algebras in the complex domain. \newline
\textbf{Applications:}\newline
\textbf{1.} In this example, we show how Noether-like operators are used to
construct invariants for the systems. We start with the system of two
second-order non-linear ODEs%
\begin{eqnarray}
ff^{\prime \prime }-gg^{\prime \prime } =e^{-f^{\prime }}\cos (g^{\prime }),
\nonumber \\
fg^{\prime \prime }+gf^{\prime \prime } =-e^{-f^{\prime }}\sin (g^{\prime }). \label{eq17}
\end{eqnarray}%
The above equations are equivalent to a system of EL equations (\ref{eq4}) on
employing the Lagrangians%
\begin{eqnarray}
L_{1} =e^{f^{\prime }}\cos (g^{\prime })+\frac{1}{2}\ln (f^{2}+g^{2}),  \quad
L_{2} =e^{f^{\prime }}\sin (g^{\prime })+\arctan \left (\frac{g}{f}\right ). \label{eq18}
\end{eqnarray}%
Now the Noether-like symmetry conditions (\ref{eq6}) for the system of ODEs has
the form%
\begin{eqnarray}
A_{1x}+f^{\prime }(A_{1f}+A_{2g})-g^{\prime }(A_{2f}-A_{1g})
=\frac{\chi
_{1f}+\chi _{2g}}{f^{2}+g^{2}}+  \nonumber \\
e^{f^{\prime }}[\{\chi _{1x}+{\small \varsigma _{1x}}+f^{\prime
}(\chi _{1f}+\chi _{2g}-{\small \varsigma _{1x}+\varsigma
_{1f}+\varsigma _{2g}})-
\nonumber \\
g^{\prime }(\chi _{2f}-\chi _{1g}-{\small \varsigma
_{2x}+\varsigma _{2f}-\varsigma _{1g}})-{\small (f^{\prime
2}-g^{\prime 2})(\varsigma
_{1f}+\varsigma _{2g})+}  \nonumber \\
{\small 2f^{\prime }g^{\prime }(\varsigma _{2f}-\varsigma
_{1g})\}\cos (g^{\prime })-}\{\chi _{2x}+{\small \varsigma
_{2x}}+g^{\prime }(\chi
_{1f}+\chi _{2g}-  \nonumber \\
{\small \varsigma _{1x}+\varsigma _{1f}+\varsigma
_{2g}})+f^{\prime }(\chi _{2f}-\chi _{1g}-{\small \varsigma
_{2x}+\varsigma _{2f}-\varsigma _{1g})-}
\nonumber \\
{\small 2f^{\prime }g^{\prime }(\varsigma _{2f}-\varsigma _{1g})}%
-(f^{\prime 2}-g^{\prime 2})({\small \varsigma _{1f}+\varsigma
_{2g})\}\sin(
g^{\prime })]+}  \nonumber \\
\frac{1}{2}\ln (f^{2}+g^{2})[{\small \varsigma _{1x}+f^{\prime
}(\varsigma _{1f}+\varsigma _{2g})-g^{\prime }(\varsigma
_{2f}-\varsigma
_{1g})]-}  \nonumber \\
\arctan (g/f)[{\small \varsigma _{2x}+g^{\prime }(\varsigma
_{1f}+\varsigma _{2g})+f^{\prime }(\varsigma _{2f}-\varsigma
_{1g})],}
\nonumber \\
\label{eq19} \\
A_{2x}+f^{\prime }(A_{2f}-A_{1g})+g^{\prime }(A_{1f}+A_{2g})
=\frac{\chi
_{2f}-\chi _{1g}}{f^{2}+g^{2}}+  \nonumber \\
e^{f^{\prime }}[\{\chi _{1x}+{\small \varsigma _{1x}}+f^{\prime
}(\chi _{1f}+\chi _{2g}-{\small \varsigma _{1x}+\varsigma
_{1f}+\varsigma _{2g}})-
\nonumber \\
g^{\prime }(\chi _{2f}-\chi _{1g}-{\small \varsigma
_{2x}+\varsigma _{2f}-\varsigma _{1g}})-{\small (f^{\prime
2}-g^{\prime 2})(\varsigma
_{1f}+\varsigma _{2g})+}  \nonumber \\
{\small 2f^{\prime }g^{\prime }(\varsigma _{2f}-\varsigma _{1g})\}\sin (g}%
^{\prime })+\{\chi _{2x}+{\small \varsigma _{2x}}+g^{\prime }(\chi
_{1f}+\chi
_{2g}-  \nonumber \\
{\small \varsigma _{1x}+\varsigma _{1f}+\varsigma
_{2g}})+f^{\prime }(\chi _{2f}-\chi _{1g}-{\small \varsigma
_{2x}+\varsigma _{2f}-\varsigma _{1g})-}
\nonumber \\
{\small 2f^{\prime }g^{\prime }(\varsigma _{2f}-\varsigma _{1g})}%
-(f^{\prime 2}-g^{\prime 2})({\small \varsigma _{1f}+\varsigma _{2g})\}\cos (g%
}^{\prime })]+  \nonumber \\
\arctan (g/f)[{\small \varsigma _{1x}+f^{\prime }(\varsigma
_{1f}+\varsigma _{2g})-g^{\prime }(\varsigma _{2f}-\varsigma
}_{1g})]+
\nonumber \\
\frac{1}{2}\ln (f^{2}+g^{2})[{\small \varsigma _{2x}+g^{\prime
}(\varsigma _{1f}+\varsigma _{2g})+f^{\prime }(\varsigma
_{2f}-\varsigma
_{1g})].}   \label{eq20}
\end{eqnarray}%
On comparing the coefficients of all independent quantities in the
above equations, we get the following system of linear PDEs
\begin{eqnarray}
{\small \varsigma _{1f}}+{\small \varsigma }_{2g} =0, \quad 
{\small \varsigma _{2f}-\varsigma }_{1g} =0, \label{eq21} \\
\chi _{1x}+{\small \varsigma }_{1x} =0,  \quad 
\chi _{2x}+{\small \varsigma _{2x}} =0, \label{eq22}\\
\chi _{1f}+\chi _{2g}-{\small \varsigma _{1x}} =0, \label{eq23}\\
\chi _{2f}-\chi _{1g}-{\small \varsigma _{2x}} =0, \label{eq24}\\
A_{1f}+A_{2g} =0,  \quad
A_{2f}-A_{1g} =0, \label{eq25}\\
\frac{f\chi _{1}+g\chi _{2}}{f^{2}+g^{2}}+\frac{1}{2}%
\ln (f^{2}+g^{2}){\small \varsigma }_{1x}-\arctan (g/f){\small \varsigma }_{2x}
=A_{1x},
\nonumber \\
\frac{f\chi _{2}-g\chi _{1}}{f^{2}+g^{2}}+\frac{1}{2}
\ln (f^{2}+g^{2}){\small \varsigma }%
_{2x}+\arctan (g/f){\small \varsigma }_{1x}
=A_{2x}. \label{eq26}
\end{eqnarray}%
Eqs. (\ref{eq21}) and (\ref{eq25}) upon using the 
analyticity of $\varsigma _{1},~\varsigma _{2},~A_1$ and $A_2$ imply that%
\begin{eqnarray}
{\small \varsigma _{1}} =  \varsigma_{1}(x), \quad {\small \varsigma _{2}} = \varsigma_{2}(x), \quad
A_1=A_1(x),\quad A_2=A_2(x).  \nonumber
\end{eqnarray}%
Eqs. (\ref{eq23}) and (\ref{eq24}) on utilizing above yield%
\begin{eqnarray}
\chi _{1} =f \varsigma_{1}^{\prime }-g \varsigma_{2}^{\prime }+g_{1}(x), \quad
\chi _{2} =f \varsigma_{2}^{\prime }+g \varsigma_{1}^{\prime }+g_{2}(x). \label{eq27}
\end{eqnarray}%
Insertion of (\ref{eq27}) in (\ref{eq22}) gives
\begin{eqnarray}
\chi _{1} =(f-x)C_{1}-gC_{2},  \quad
\chi _{2} =(f-x)C_{2}+gC_{1}. \label{eq28}
\end{eqnarray}%
Now by differentiating (\ref{eq26}) with respect of $f$ and 
$g$ and using Eqs. (\ref{eq28}) in it gives us the
following solution $C_{1}=\chi _{1}=0=\chi _{2}=C_{2},$ and ${\small \varsigma _{1}}=C_{3},$\ $%
{\small \varsigma _{2}}=C_{4}{\small ,}$ while the gauge functions
are determined to be $A_{1}=C_{5},$\ $A_{2}=C_{6}.$ Therefore, we
obtain a single Noether-like operator which is translation in $x$ of
system (\ref{eq17}). Hence, by utilizing (\ref{eq12}) the invariants of (\ref{eq17})
are
\begin{eqnarray}
I_{1} =e^{f^{\prime }}\cos (g^{\prime })+\frac{1}{2}\ln
(f^{2}+g^{2})-e^{f^{\prime }}(f^{\prime }\cos (g^{\prime })-g^{\prime }\sin
(g^{\prime })),  \nonumber \\
I_{2} =e^{f^{\prime }}\sin (g^{\prime })+\arctan (g/f)-e^{f^{\prime
}}(f^{\prime }\sin (g^{\prime })+g^{\prime }\cos (g^{\prime })). \label{eq29}
\end{eqnarray}%
This case corresponds to the one-dimensional realization of the complex algebra
spanned by $\partial_{x}$ mentioned above. In fact system (\ref{eq17}) can be converted into an r-CODE
$uu''=e^{-u'}$, which has Lagrangian $L=e^{u'}+\log u$.\newline
\textbf{2.} In this example, we study the invariant properties of the system
\begin{eqnarray}
2f^{\prime \prime }+f^{\prime }f^{\prime \prime }-g^{\prime }g^{\prime
\prime } =xe^{-f^{\prime }}\cos (g^{\prime }),  \nonumber \\
2g^{\prime \prime }+f^{\prime }g^{\prime \prime }+g^{\prime }f^{\prime
\prime } =-xe^{-f^{\prime }}\sin (g^{\prime }). \label{eq30}
\end{eqnarray}%
It can be checked that the operators $\partial_{f}$ and $\partial_{g}$ 
satisfy (\ref{eq6}) which thus provide two Noether-like operators for the
above system relative to
\begin{eqnarray}
L_{1} =e^{f^{\prime} }(f^{\prime }\cos (g^{\prime })-g^{\prime }\sin (g^{\prime}))+xf,  \nonumber \\
L_{2} =e^{f^{\prime} }(f^{\prime }\sin (g^{\prime })+g^{\prime }\cos (g^{\prime}))+xg, \label{eq31}
\end{eqnarray}%
by solving Noether-like symmetry conditions (\ref{eq6})%
\begin{eqnarray}
A_{1x}+f^{\prime }(A_{1f}+A_{2g})-g^{\prime }(A_{2f}-A_{1g})=f\small \varsigma  _{1}-g\small \varsigma
_{2}+\chi _{1}x+[\{\chi _{1x}+   \nonumber \\
f^{\prime }(\chi _{1f}+\chi _{2g})-g^{\prime
}(\chi _{2f}-\chi _{1g})-f^{\prime }(\small \varsigma  _{1f}+\small \varsigma  _{2g})+g^{\prime }(\small \varsigma  _{2f}-\small \varsigma
_{1g})-  \nonumber \\
(f^{\prime 2}-g^{\prime 2})(\small \varsigma  _{1f}+\small \varsigma  _{2g})+2f^{\prime
}g^{\prime }(\small \varsigma  _{2f}-\small \varsigma  _{1g})\}e^{f^{\prime }}(\cos g^{\prime }(1+f^{\prime })-g^{\prime }\sin g^{\prime })
\nonumber \\
-\{\chi _{2x}+g^{\prime }(\chi _{1f}+\chi _{2g})+f^{\prime }(\chi _{2f}-\chi
_{1g})+f^{\prime }(\small \varsigma  _{2f}-\small \varsigma  _{1g})-g^{\prime }(\small \varsigma  _{1f}+\small \varsigma
_{2g})-  \nonumber \\
2f^{\prime }g^{\prime }(\small \varsigma  _{1f}-\small \varsigma  _{2g})-(f^{\prime 2}-g^{\prime 2})(\small \varsigma  _{2f}-\small \varsigma  _{1g})\}e^{f^{\prime
}}((1+gf^{\prime })\sin g^{\prime }+g^{\prime }\cos g^{\prime })]+  \nonumber \\
x(f\small \varsigma _{1x}-g\small \varsigma  _{2x})+x(ff^{\prime }-gg^{\prime })(\small \varsigma  _{1f}+\small \varsigma  _{2g})-x(fg^{\prime }+f^{\prime
}g)(\small \varsigma  _{2f}-\small \varsigma  _{1g})  \nonumber \\
e^{f^{\prime }}[(f^{\prime }\cos g^{\prime }-g^{\prime }\sin g^{\prime })\small \varsigma
_{1x}-(f^{\prime }\sin g^{\prime }+g^{\prime }\cos g^{\prime })\small \varsigma  _{2x} +((f^{\prime 2}-g^{\prime 2})\cos g^{\prime }-
\nonumber \\
2f^{\prime }g^{\prime }\sin
g^{\prime })(\small \varsigma  _{1f}+\small \varsigma  _{2g})-
((f^{\prime 2}-g^{\prime 2})\sin g^{\prime
}+2f^{\prime }g^{\prime }\cos g^{\prime })(\small \varsigma  _{2f}-\small \varsigma  _{1g})]. \nonumber 
\\ \label{eq32}\\
A_{2x}+g^{\prime }(A_{1f}+A_{2g})+f^{\prime }(A_{2f}-A_{1g})=g\small \varsigma  _{1}+f\small \varsigma
_{2}+\chi _{2}x+[\{\chi _{1x}+  \nonumber \\
f^{\prime }(\chi _{1f}+\chi _{2g})-g^{\prime
}(\chi _{2f}-\chi _{1g})-f^{\prime }(\small \varsigma  _{1f}+\small \varsigma  _{2g})+g^{\prime }(\small \varsigma  _{2f}-\small \varsigma
_{1g})-  \nonumber \\
(f^{\prime 2}-g^{\prime 2})(\small \varsigma  _{1f}+\small \varsigma  _{2g})+2f^{\prime
}g^{\prime }(\small \varsigma  _{2f}-\small \varsigma  _{1g})\}e^{f^{\prime }}((1+f^{\prime })\sin g^{\prime }+g^{\prime }\cos g^{\prime })
\nonumber \\
+\{\chi _{2x}+g^{\prime }(\chi _{1f}+\chi _{2g})+f^{\prime }(\chi _{2f}-\chi
_{1g})+f^{\prime }(\small \varsigma  _{2f}-\small \varsigma  _{1g})-g^{\prime }(\small \varsigma  _{1f}+\small \varsigma
_{2g})-  \nonumber \\
2f^{\prime }g^{\prime }(\small \varsigma  _{1f}+\small \varsigma  _{2g})-(f^{\prime 2}-g^{\prime 2})(\small \varsigma  _{2f}-\small \varsigma  _{1g})\}e^{f^{\prime
}}((1+f^{\prime })\cos g^{\prime }-g^{\prime }\sin g^{\prime })]+  \nonumber \\
x(g\small \varsigma _{1x}+f\small \varsigma  _{2x})+x(ff^{\prime }-gg^{\prime })(\small \varsigma  _{2f}-\small \varsigma  _{1g})+x(fg^{\prime }+f^{\prime
}g)(\small \varsigma  _{1f}+\small \varsigma  _{2g})  \nonumber \\
e^{f^{\prime }}[(f^{\prime }\cos g^{\prime }-g^{\prime }\sin g^{\prime })\small \varsigma
_{2x}+(f^{\prime }\sin g^{\prime }+g^{\prime }\cos g^{\prime })\small \varsigma  _{1x}
\nonumber \\
+(f^{\prime 2}-g^{\prime 2})(\sin g^{\prime }(\small \varsigma  _{1f}+\small \varsigma  _{2g})+\cos
g^{\prime }(\small \varsigma  _{2f}-\small \varsigma  _{1g}))  \nonumber \\
+2f^{\prime }g^{\prime }(\cos g^{\prime }(\small \varsigma  _{1f}+\small \varsigma  _{2g})-\sin
g^{\prime }(\small \varsigma  _{2f}-\small \varsigma  _{1g})].\nonumber \\\label{eq33}
\end{eqnarray}%
Again, by comparison of the coefficients of independent functions, the above equations give rise to
\begin{eqnarray}
\small \varsigma  _{1x}=0,\quad \small \varsigma  _{2x}=0,    \label{eq34}\\ 
\small \varsigma  _{1f}+\small \varsigma  _{2g}=0,\quad \small \varsigma  _{2f}-\small \varsigma  _{1g}=0,   \label{eq35}\\ 
\chi _{1x}=0,\quad \chi _{2x}=0, \label{eq36}\\ 
\chi _{1f}+\chi _{2g}=0,\quad \chi _{2f}-\chi _{1g}=0,   \label{eq37}\\ 
\small \varsigma  _{1}f-\small \varsigma  _{2}g+\chi _{1}x=A_{1x},  \label{eq38}\\ 
\small \varsigma  _{1}g+\small \varsigma  _{2}f+\chi _{2}x=A_{2x},   \label{eq39}\\ 
A_{1f}+A_{2g}=0,\quad A_{2f}-A_{1g}=0. \label{eq40}
\end{eqnarray}%
Equations (\ref{eq34})$-$(\ref{eq37}) yield $\small \varsigma  _{1}=C_{1},$ $\small
\varsigma  _{2}=C_{2}$ and $\chi _{1}=C_{3},$ $\chi _{2}=C_{4}.$
From equations (\ref{eq38}) and (\ref{eq39}), we obtain $\small \varsigma  _{1}=\small
\varsigma  _{2}=0.$ Thus in this case we have $\partial_{f}$
and $\partial_{g}$ relative to the gauge functions
$A_{1}=x^{2}/2,$ while $A_{2}=0$, respectively. Invoking
(\ref{eq12}) the two first integrals
\begin{eqnarray}
I_{1} =e^{f^{\prime }}\cos (g^{\prime })+e^{f^{\prime }}(f^{\prime }\cos
g^{\prime }-g^{\prime }\sin (g^{\prime }))-x^{2}/2,  \nonumber \\
I_{2} =e^{f^{\prime }}\sin (g^{\prime })+e^{f^{\prime }}(f^{\prime }\sin
g^{\prime }+g^{\prime }\cos (g^{\prime })), \label{eq41}
\end{eqnarray}%
of (\ref{eq30}) can be deduced. Notice that the two operators $\partial_{f}$ 
and $\partial_{g}$ correspond to a single complex symmetry $%
\partial_{u}$. Also system (\ref{eq30}) arises from the EL equation $(2+u')u''=xe^{-u'},$ with Lagrangian
$L=u'e^{u'}+xu$. \newline
\textbf{3.} Here we arrive at three Noether-like operators of the following
systems of ODEs%
\begin{eqnarray}
(2+4f^{\prime }+f^{\prime 2}-g^{\prime 2})f^{\prime \prime }-2(2g^{\prime
}+f^{\prime }g^{\prime })g^{\prime \prime } =e^{-f^{\prime }}\cos
g^{\prime },  \nonumber \\
(2+4f^{\prime }+f^{\prime 2}-g^{\prime 2})g^{\prime \prime }+2(2g^{\prime
}+f^{\prime }g^{\prime })f^{\prime \prime } =-e^{-f^{\prime }}\sin
g^{\prime }, \label{eq42}
\end{eqnarray}%
that has a variational structure and admits Lagrangians%
\begin{eqnarray}
L_{1} =(f^{\prime 2}-g^{\prime 2})e^{f^{\prime }}\cos g^{\prime
}-2f^{\prime }g^{\prime }e^{f^{\prime }}\sin (g^{\prime })+f,  \nonumber \\
L_{2} =2f^{\prime }g^{\prime }e^{f^{\prime }}\cos (g^{\prime })+(f^{\prime
2}-g^{\prime 2})e^{f^{\prime }}\sin (g^{\prime })+g. \label{eq43}
\end{eqnarray}%
The system (\ref{eq42}) admits the Noether-like operators $\partial_{x}$, 
$\partial_{f}$ and $\partial_{g}$
corresponding to the above Lagrangians. By using conditions (\ref{eq12})
the four first integrals are
found to be%
\begin{eqnarray}
I_{1} =f-e^{f^{\prime }}((f^{\prime 2}-g^{\prime
2})\cos g^{\prime }-2f^{\prime }g^{\prime }\sin g^{\prime })-e^{f^{\prime }}
((f^{\prime 3}-3f^{\prime }g^{\prime 2})\cos g^{\prime
}-  \nonumber \\
(3f^{\prime 2}g^{\prime }-g^{\prime 3}))\sin g^{\prime },  \nonumber \\
I_{2} =g-e^{f^{\prime }}((f^{\prime 2}-g^{\prime
2})\sin g^{\prime }+2f^{\prime }g^{\prime }\cos g^{\prime })-e^{f^{\prime }}((f^{\prime 3}-3f^{\prime }g^{\prime 2})\sin g^{\prime
}+  \nonumber \\
(3f^{\prime 2}g^{\prime }-g^{\prime 3}))\cos g^{\prime },  \nonumber \\
I_{3} =2e^{f^{\prime }}(f^{\prime }\cos g^{\prime }-g^{\prime }\sin
g^{\prime })+e^{f^{\prime }}((f^{\prime 2}-g^{\prime 2})\cos g^{\prime
}-2f^{\prime }g^{\prime }\sin g^{\prime })-x,  \nonumber \\
I_{4} =2e^{f^{\prime }}(f^{\prime }\sin g^{\prime }+g^{\prime }\cos
g^{\prime })+e^{f^{\prime }}((f^{\prime 2}-g^{\prime 2})\sin g^{\prime
}+2f^{\prime }g^{\prime }\cos g^{\prime }).\nonumber \\ \label{eq44}
\end{eqnarray}%
of (\ref{eq42}). Here, the first two invariants correspond to the
Noether-like operator $\partial /\partial x$. In this case three
operators can be transformed into a realization of the
$2-$dimensional Abelian complex algebra. The r-CODE here is
$(2+4u'+u'^2)u''=e^{-u'},$ which has Lagrangian $L=u'^2e^{u'}+u$.
\newline
\textbf{4.} Now we look for a system which respects the second
two-dimensional algebra in the complex domain. It  helps us in
constructing four Noether-like operators. For this we consider
\begin{eqnarray}
x(2+f^{\prime })f^{\prime \prime }-xg^{\prime }g^{\prime \prime }
=f^{\prime }+1,  \nonumber \\
x(2+f^{\prime })g^{\prime \prime }+xg^{\prime }f^{\prime \prime }=g^{\prime }, \label{eq45}
\end{eqnarray}%
which admits four Noether-like operators
\begin{eqnarray}
\mathbf{X}_{1}=\frac{\partial }{\partial f},\quad \mathbf{X}_{2}=\frac{\partial }{%
\partial g},\quad \mathbf{X}_{3}=x\frac{\partial }{\partial x}+f\frac{\partial }{%
\partial f}+g\frac{\partial }{\partial g},\nonumber \\
\mathbf{X}_{4}=g\frac{\partial }{%
\partial f}-f\frac{\partial }{\partial g}, \label{eq46}
\end{eqnarray}%
with respect to the Lagrangians
\begin{eqnarray}
L_{1} =\frac{e^{f^{\prime }}}{x}(f^{\prime }\cos (g^{\prime })-g^{\prime
}\sin (g^{\prime })),  \quad
L_{2} =\frac{e^{f^{\prime }}}{x}(g^{\prime }\cos (g^{\prime })+f^{\prime
}\sin (g^{\prime })). \label{eq47}
\end{eqnarray}%
By utilizing (\ref{eq12}) we determine the following four real first integrals
\begin{eqnarray}
I_{1} =\left((1+f^{\prime })\cos (g^{\prime
})-g^{\prime }\sin (g^{\prime })\right)e^{f^{\prime }}/x,  \nonumber\\
I_{2} =\left((1+f^{\prime })\sin (g^{\prime
})+g^{\prime }\cos (g^{\prime })\right)e^{f^{\prime }}/x,  \nonumber \\
I_{3} =((f(1+f^{\prime
})-gg^{\prime }-x(f^{\prime 2}-g^{\prime 2}))\cos (g^{\prime })- (g(1+f^{\prime })+ \nonumber \\
fg^{\prime }-2xf^{\prime }g^{\prime
})\sin (g^{\prime }))e^{f^{\prime }}/x,  \nonumber \\
I_{4} =((f(1+f^{\prime
})-gg^{\prime }-x(f^{\prime 2}-g^{\prime 2}))\sin (g^{\prime })+ (g(1+f^{\prime })+ \nonumber \\
fg^{\prime }-2xf^{\prime }g^{\prime
})\cos (g^{\prime }))e^{f^{\prime }}/x. \nonumber \\ \label{eq48}
\end{eqnarray}%
for the system (\ref{eq45}).

The above four Noether-like operators are also Noether symmetries of
the system (\ref{eq45}) as these satisfy (\ref{eq10}) and (\ref{eq11}). It would have
been more difficult to determine these via the usual Noether approach especially $%
\mathbf{X}_{3}$ and $\mathbf{X}_{4},$ which are mainly the result of complex
encoding. This enters us into an open domain of nice and interesting
problems in the complex domain with the clear indication that the solution
of each or any of these would result in something remarkable in the real
domain. Therefore, it is indispensable to classify such systems of ODEs with
respect to the Noether-like operators they admit. The r-CODE is $x(2+u')u''=1+u',$ which has Lagrangian
$L=u'e^{u'}/x$.

In \cite{18}, it is shown that the $5-$dimensional Noether algebra of a free
particle equation is a subalgebra of $sl(3,\Re)$ that give rise to five invariants. The Lagrangian that yields the maximum
number of Noether symmetries is known as the standard Lagrangian. The concept of
an alternative Lagrangian could be exercised to construct other invariants
of differential equations. The invariants relative to such Lagrangians are
expressed in terms of invariants of the standard Lagrangian. A similar situation
arises in the case of systems of two ODEs. In the next example, we highlight
the significance of an alternative Lagrangian in the variational problem of
systems of ODEs. \newline
\textbf{Alternative Lagrangians:}\newline
The use of Lagrangians $L_{1}=e^{f^{\prime }}\cos (g^{\prime })+f,$ and $%
L_{2}=e^{f^{\prime }}\sin (g^{\prime })+g,$ in the EL-equations $(4)$
yields the system
\begin{eqnarray}
f^{\prime \prime } =e^{-f^{\prime }}\cos (g^{\prime }),  \quad
g^{\prime \prime } =-e^{-f^{\prime }}\sin (g^{\prime }). \label{eq49}
\end{eqnarray}%
These Lagrangians admit three Noether-like operators
\begin{equation}
\mathbf{X}_{1}=\frac{\partial }{\partial x},\quad \mathbf{X}_{2}=\frac{\partial }{%
\partial f},\quad \mathbf{X}_{3}=\frac{\partial }{\partial g}, \label{eq50}
\end{equation}%
which provide the first integrals
\begin{eqnarray}
I_{1} =x-e^{f^{\prime }}\cos (g^{\prime }),  \quad
I_{2} =e^{f^{\prime }}\sin (g^{\prime }),  \nonumber \\
I_{3} =e^{f^{\prime }}\left(f^{\prime }\cos (g^{\prime })-g^{\prime
}\sin (g^{\prime })-\cos (g^{\prime })\right)-f,  \nonumber \\
I_{4} =e^{f^{\prime }}\left(g^{\prime }\cos (g^{\prime })+f^{\prime
}\sin (g^{\prime })-\sin (g^{\prime })\right)-g, \label{eq51}
\end{eqnarray}%
for (\ref{eq49}). The system (\ref{eq49}) also has an alternative
Lagrangians%
\begin{eqnarray}
\tilde{L}_{1} =-\frac{1}{2x}\ln \left(1+x^{2}e^{-2f^{\prime
}}-2xe^{-f^{\prime }}\cos (g^{\prime })\right),  \nonumber \\
\tilde{L}_{2} =\frac{-1}{x}\arctan \left (\frac{xe^{-f^{\prime }}\sin
(g^{\prime })}{1-xe^{-f^{\prime }}\cos (g^{\prime })}\right ), \label{eq52}
\end{eqnarray}
which possess the Noether-like operators
\begin{eqnarray}
\mathbf{X}_{4} =x\frac{\partial }{\partial x}+(x+f)\frac{\partial }{%
\partial f}+g\frac{\partial }{\partial g},  \quad
\mathbf{X}_{5} =g\frac{\partial }{\partial f}-(x+f)\frac{\partial }{%
\partial g}. \label{eq53}
\end{eqnarray}%
The above operators reveal two more first integrals
\begin{eqnarray}
I_{5} =(1/2)\ln [(x-e^{f^{\prime }}\cos (g^{\prime }))^{2}+e^{2f^{\prime
}}\sin ^{2}g^{\prime }]+[(x-e^{f^{\prime }}\cos (g^{\prime }))  \nonumber \\
\{e^{f^{\prime }}(f^{\prime }\cos (g^{\prime })-g^{\prime }\sin
g^{\prime })-f^{\prime }-f\}]-e^{f^{\prime }}\sin (g^{\prime })\{e^{f^{\prime }}(g^{\prime
}\cos (g^{\prime })+  \nonumber \\
f^{\prime }\sin (g^{\prime }))-g^{\prime }-g\}/((x-e^{f^{\prime }}\cos g^{\prime
})^{2}+e^{2f^{\prime }}\sin ^{2}g^{\prime }),  \nonumber \\
I_{6} =\arctan (e^{f^{\prime }}\sin (g^{\prime })/e^{f^{\prime }}\cos
g^{\prime }-x)+[(x-e^{f^{\prime }}\cos (g^{\prime }))  \nonumber \\
\{e^{f^{\prime }}(g^{\prime
}\cos (g^{\prime })+
f^{\prime }\sin (g^{\prime }))-g^{\prime }-g\}]+e^{f^{\prime }}
\sin (g^{\prime })\{e^{f^{\prime }}(f^{\prime }\cos (g^{\prime })\nonumber\\
-g^{\prime }\sin
g^{\prime })-f^{\prime }-f\}/((x-e^{f^{\prime }}\cos g^{\prime
})^{2}+e^{2f^{\prime }}\sin ^{2}g^{\prime }), \label{eq54}
\end{eqnarray}%
for the system (\ref{eq49}). It is important to see that alternative
Lagrangians can also be used to reveal other first integrals of DEs.
Perhaps more importantly from a physical point of view, alternative
Lagrangians allow us to have a complete determination of all the
physical constants or integrals (see \cite{27}).
\section{Classification of Noether-like operators}
We firstly present the classification of three Noether point
symmetries for complex Lagrangians on the line. Table 1 lists
different cases of $3-$dimensional Noether complex algebras together
with their Lagrangians and representative EL r-CODEs. These are
taken from \cite{19}. 

We then classify systems of ODEs with respect to $5$ and $6$
Noether-like operators using CAS. We divide our discussion into two parts.
There are three cases in which the commutators of Noether-like
operators are closed which are described in the first part. The
other part contains the remaining cases. We obtain the
classification of $5-$ and $6-$dimensional Noether-like operators.
These cases correspond to the cases described in Table 1. Moreover,
in each case we write down the first integrals associated with the
Noether-like operators which reveal an important feature of such
operators.

Our purpose is to find out those systems of ODEs that can be mapped to
r-CODEs via complex transformations. In particular, if a system of two ODEs
with Lagrangians admitting at least $5$ Noether-like operators
is transformable to an r-CODE then the symmetry structure of that r-CODE
will correspond to one of the cases described in Table 1. For instance,
notice that the case $\mathcal{N}_{3,5}^{1}$ correspond to the first case in
Table $1$ as the Noether-like operators $\mathbf{X}_{1},...,\mathbf{X}_{5}
$ belong to complex symmetries $\partial _{x},\partial _{u}$ and $x\partial
_{x}-u\partial _{u}.$ Similarly, all cases described hereafter can
be projected to each case in Table $1$. This will help us in
understanding the inverse problem for systems of two ODEs via r-CODEs.
\newline\newline
Table 1. Classification of Noether Algebras for Lagrangians on the Line%
\newline
\begin{tabular}{|l|l|l|l|}
\hline
No.  ~~Noether symmetries ~~Lagrangians  ~~~~~~~~~~~~~~~~ Representative ELEs \\ \hline
$\mathcal{N}_{3,5}^{1}$
\begin{tabular}{l}
$\mathbf{Z}_{1}=\partial _{x},\mathbf{Z}_{2}=\partial _{u},$ \\
$\mathbf{Z}_{3}=x\partial _{x}-u\partial _{u}$%
\end{tabular}
~~~~$L=-4u^{\prime 1/2}+u$ ~~~~~~~~~~~ $u^{\prime \prime }=u^{\prime 3/2}$ \\ \hline
$\mathcal{N}_{3,5}^{2}$
\begin{tabular}{l}
$\mathbf{Z}_{1}=\partial _{x},\mathbf{Z}_{2}=\partial _{u},$ \\
$\mathbf{Z}_{3}=u\partial _{x}-x\partial _{u}$%
\end{tabular}
~~~~$L=-(1+u^{\prime 2})^{1/2}+xu^{\prime }$ ~~$u^{\prime \prime
}=(1+u^{\prime 2})^{3/2}$ \\ \hline
$\mathcal{N}_{3,6}^{1}$
\begin{tabular}{l}
$\mathbf{Z}_{1}=\partial _{u},$ \\
$\mathbf{Z}_{2}=x\partial _{x}+u\partial _{u}$ \\
$\mathbf{Z}_{3}=2xu\partial _{x}+u^{2}\partial _{u}$%
\end{tabular}
~~$L=\frac{u^{\prime }}{x}+\frac{1}{2xu^{\prime }}$ ~~~~~~~~~~~~~~ $xu^{\prime \prime
}=u^{\prime 3}-\frac{1}{2}u^{\prime }$ \\ \hline
$\mathcal{N}_{3,5}^{3}$
\begin{tabular}{l}
$\mathbf{Z}_{1}=\partial _{x},\mathbf{Z}_{2}=\partial _{u},$ \\
$\mathbf{Z}_{3}=(bx+u)\partial _{x}+$ \\
$(bu-x)\partial _{u}$%
\end{tabular}
\begin{tabular}{l}
~~$L=e^{-b\arctan u^{\prime }}/(1+b^{2})$ \\
$(1+u^{\prime 2})^{1/2}+xu^{\prime }$%
\end{tabular}
~$u^{\prime \prime }=(1+u^{\prime 2})^{\frac{3}{2}}e^{b\arctan u^{\prime }}$
\\ \hline
$\mathcal{N}_{3,6}^{2}$
\begin{tabular}{l}
$\mathbf{Z}_{1}=\partial _{u},$ \\
$\mathbf{Z}_{2}=x\partial _{x}+u\partial _{u}$ \\
$\mathbf{Z}_{3}=2xu\partial _{x}+$ \\
$(u^{2}-x^{2})\partial _{u}$%
\end{tabular}
 ~~~~~~$L=\frac{\sqrt{1+u^{\prime 2}}}{x}+\frac{Au^{\prime }}{x}$
\begin{tabular}{l}
~~~~~~~~~$xu^{\prime \prime }=u^{\prime 3}+u^{\prime }+$ \\
~~~~~~~~~$A(1+u^{\prime 2})^{3/2}$%
\end{tabular}
\\ \hline
\end{tabular}%
\newline
\newline
We now elaborate on the naming used. The indices in $\mathcal{N}%
_{i,j}^{a},$ are set as follows: $a$ corresponds to different cases, $i$
stands for the dimension of the complex algebras and $j$ denotes the type of
Noether-like operators that can arise from these cases. It may be noticed
that we obtain four cases of $3-$dimensional complex algebra that give rise
to 5$-$Noether-like operators corresponding to $\mathcal{N}_{3,5}^{1},~
\mathcal{N}_{3,5}^{2},~\mathcal{N}_{3,5}^{3}.$ Similarly, $\mathcal{N}%
_{3,6}^{1}$ and $\mathcal{N}_{3,6}^{2}$ are  cases of $6-$%
Noether-like operators. There is one more case $\mathcal{N}_{3,6}^{3}$ of $%
6-$Noether-like operators corresponding to $3-$dimensional complex
Noether algebra. Due to its length, we have summarized it at the end
of this section.\newline
\newline
\textbf{Closed Noether-like Operators:}\newline
\textbf{Case }$\mathcal{N}_{3,5}^{1}$ \newline
We first consider the system of nonlinear ODEs
\begin{eqnarray}
f^{\prime \prime } =(f^{\prime 2}+g^{\prime 2})^{3/4}\cos \left(\frac{3}{2}%
\arctan  \left (\frac{g^{\prime }}{f^{\prime }}\right)\right),  \nonumber \\
g^{\prime \prime } =(f^{\prime 2}+g^{\prime 2})^{3/4}\sin \left(\frac{3}{2}%
\arctan  \left (\frac{g^{\prime }}{f^{\prime }}\right)\right),\label{eq55}
\end{eqnarray}%
which can be obtained from the Lagrangians%
\begin{eqnarray}
L_{1} =-4(f^{\prime 2}+g^{\prime 2})^{1/4}\cos \left(\frac{1}{2} \arctan \left( 
\frac{g^{\prime}}{f^{\prime}}\right) \right) +f,  \nonumber \\
L_{2} =-4(f^{\prime 2}+g^{\prime 2})^{1/4}\sin \left(\frac{1}{2} \arctan \left( 
\frac{g^{\prime}}{f^{\prime}}\right) \right) +g,\label{eq56}
\end{eqnarray}%
with the aid of (\ref{eq4}). This system corresponds to the r-CODE listed for this case in Table $1$.
 The system (\ref{eq55}) admits the following five
Noether-like operators
\begin{eqnarray}
\mathbf{X}_{1} =\frac{\partial }{\partial x},\quad \mathbf{X}_{2}=\frac{\partial
}{\partial f},\quad \mathbf{X}_{3}=\frac{\partial }{\partial g},\quad \mathbf{X}_{5}=-g%
\frac{\partial }{\partial f}+f\frac{\partial }{\partial g},  \nonumber \\
\mathbf{X}_{4} =x\frac{\partial }{\partial x}-f\frac{\partial }{\partial f%
}-g\frac{\partial }{\partial g},\label{eq57}
\end{eqnarray}%
corresponding to the above Lagrangians which can be verified from (\ref{eq6}).
 Since the first three operators are
translations in $x,f$ and $g$ they therefore leave both the system and the
Lagrangians invariant. Consequently, these turn out to be Noether
symmetries. The use of operator $\mathbf{X}_{1}$ in (\ref{eq12}) gives
\begin{eqnarray}
I_{1} =2r^{1/2}\cos (\theta /2)-f,  \quad
I_{2} =2r^{1/2}\sin (\theta /2)-g, \label{eq58}
\end{eqnarray}%
i.e., the two first integrals of the system (\ref{eq55}) where $r$ and
$\theta $ are defined as%
\begin{eqnarray}
r^{2} =f^{\prime 2}+g^{\prime 2},  \quad
\theta  =\arctan (g^{\prime }/f^{\prime }).\label{eq59}
\end{eqnarray}%
In a similar way, $\partial /\partial f$ and $\partial /\partial g$ give the
two invariants%
\begin{eqnarray}
I_{3} =2r^{-1/2}\cos (\theta /2)-x,  \quad
I_{4} =2r^{-1/2}\sin (\theta /2), \label{eq60}
\end{eqnarray}%
for the system (\ref{eq55}). Similarly, for (\ref{eq55}) one can determine the
first integrals
\begin{eqnarray}
I_{5} =2r^{-1/2}(f\cos (\theta /2)+g\sin (\theta /2))+xf-2xr^{1/2}\cos (\theta
/2)-4,  \nonumber \\
I_{6} =2r^{-1/2}(g\cos (\theta /2)-f\sin (\theta /2))+xg-2xr^{1/2}\sin (\theta
/2),\label{eq61}
\end{eqnarray}%
corresponding to $\mathbf{X}_{4}$ and $\mathbf{X}_{5}$. We claim
that any system exhibiting these operators can be mapped to the
system (\ref{eq55}). The commutators of these operators are
\begin{eqnarray}
\lbrack \mathbf{X}_{1},\mathbf{X}_{2}] =[\mathbf{X}_{2},\mathbf{X}_{3}]=[%
\mathbf{X}_{1},\mathbf{X}_{3}]=[\mathbf{X}_{1},\mathbf{X}_{5}]=[\mathbf{X}%
_{4},\mathbf{X}_{5}]=0,  \nonumber \\
\lbrack \mathbf{X}_{1},\mathbf{X}_{4}] =\mathbf{X}_{1},\quad [\mathbf{X}_{2},%
\mathbf{X}_{5}]=\mathbf{X}_{3},\quad [\mathbf{X}_{3},\mathbf{X}_{5}]=-\mathbf{X}%
_{2},  \nonumber \\
\lbrack \mathbf{X}_{2},\mathbf{X}_{4}] =-\mathbf{X}_{2},\quad [\mathbf{X}_{3},%
\mathbf{X}_{4}]=-\mathbf{X}_{3}. \label{eq62}
\end{eqnarray}%
Therefore, these form a closed algebra. Noteworthy, the system (\ref{eq55}) admits
the maximal $5-$dimensional Noether-like operators.\newline
\newline
\textbf{Case }$\mathcal{N}_{3,5}^{2}$\newline
Consider the pair of two coupled Lagrangians%
\begin{eqnarray}
L_{1} =-r_{1}^{1/2}\cos (\theta _{1}/2)+xf^{\prime },  \quad
L_{2} =-r_{1}^{1/2}\sin (\theta _{1}/2)+xg^{\prime }, \label{eq63}
\end{eqnarray}%
where%
\begin{eqnarray}
r_{1}^{2} =(1+f^{\prime 2}-g^{\prime 2})^{2}+4f^{\prime 2}g^{\prime 2},
\quad	
\theta _{1} =\arctan \left ( \frac{2f^{\prime }g^{\prime }}{1+f^{\prime 2}-g^{\prime 2}} \right). \label{eq64}
\end{eqnarray}%
The system that possesses these Lagrangians is
\begin{eqnarray}
f^{\prime \prime } =r_{1}^{3/2}\cos (3\theta _{1}/2), \quad
g^{\prime \prime } =r_{1}^{3/2}\sin (3\theta _{1}/2), \label{eq65}
\end{eqnarray}%
where $r_{1}$ and $\theta _{1}$ are given by (\ref{eq64}). This system
corresponds to the r-CODE of Table 1. The five Noether-like
operators admitted by the above system are
\begin{eqnarray}
\mathbf{X}_{1} =\frac{\partial }{\partial x},\quad \mathbf{X}_{2}=\frac{\partial
}{\partial f},\quad \mathbf{X}_{3}=\frac{\partial }{\partial g},  \nonumber \\
\mathbf{X}_{4} =f\frac{\partial }{\partial x}-x\frac{\partial }{\partial f},\quad %
\mathbf{X}_{5}=g\frac{\partial }{\partial x}+x\frac{\partial }{\partial g}. \label{eq66}
\end{eqnarray}%
The use of operators $\mathbf{X}_{2}$ and $\mathbf{X}_{3}$ in (\ref{eq12}) gives
rise to
\begin{eqnarray}
I_{1} =x-r_{1}^{-1/2}\left (f^{\prime }\cos (\theta _{1}/2)+g^{\prime }\sin
(\theta _{1}/2)\right ),  \nonumber \\
I_{2} =r_{1}^{-1/2}\left(g^{\prime }\cos (\theta _{1}/2)-f^{\prime }\sin
(\theta _{1}/2)\right). \label{eq67}
\end{eqnarray}%
Now corresponding to $\partial /\partial x,$ we have the first integrals
\begin{eqnarray}
I_{3} =r_{1}^{-1/2}\cos (\theta _{1}/2)+f,  \quad
I_{4} =-r_{1}^{-1/2}\sin (\theta _{1}/2)+g, \label{eq68}
\end{eqnarray}%
for (\ref{eq65}). By invocation of (\ref{eq12}) the operators $\mathbf{X}_{4}$ and $%
\mathbf{X}_{5}$ give
\begin{eqnarray}
I_{5} =\frac{x^{2}}{2}+\frac{f^{2}-g^{2}}{2}+r_{1}^{-1/2}\left ((f-xf^{\prime
})\cos (\theta _{1}/2)+(xg^{\prime }-g)\sin (\theta _{1}/2)\right ),  \nonumber \\
I_{6} =fg+r_{1}^{-1/2}\left((g-xg^{\prime
})\cos (\theta _{1}/2)-(f-xf^{\prime })\sin (\theta _{1}/2)\right), \nonumber \\ 
\label{eq69}
\end{eqnarray}%
We now calculate the commutators of these operators
\begin{eqnarray}
\lbrack \mathbf{X}_{1},\mathbf{X}_{2}] =[\mathbf{X}_{2},\mathbf{X}_{3}]=[%
\mathbf{X}_{1},\mathbf{X}_{3}]=[\mathbf{X}_{3}\mathbf{,X}_{4}]=[\mathbf{X}%
_{2},\mathbf{X}_{5}]=0, \nonumber \\
\lbrack \mathbf{X}_{1},\mathbf{%
X}_{6}]=0, \quad [ \mathbf{X}_{1},\mathbf{X}_{4}] =-\mathbf{X}_{2},\quad [\mathbf{X}%
_{2},\mathbf{X}_{4}]=\mathbf{X}_{1},\quad [\mathbf{X}_{3},\mathbf{X}_{5}]=%
\mathbf{X}_{1},  \nonumber \\
\lbrack \mathbf{X}_{2},\mathbf{X}_{6}] =\mathbf{X}_{3},\quad[\mathbf{X}_{3},\mathbf{X}_{6}]=2\mathbf{X}%
_{2},\quad[\mathbf{X}%
_{4},\mathbf{X}_{6}]=-\mathbf{X}_{5},  \nonumber \\
\lbrack \mathbf{X}_{5},\mathbf{X}_{6}] =-\mathbf{X}_{4},\quad [\mathbf{X}%
_{1},\mathbf{X}_{5}]=\mathbf{X}_{3},\quad [\mathbf{X}%
_{4},\mathbf{X}_{5}]=f\partial_{g}+g \partial_{f}=\mathbf{X}_{6}, \label{eq70}
\end{eqnarray}%
Notice that we obtain an extra operator $\mathbf{X}_6$ which is not a Noether-like
operator. The nice feature about it is that it forms a closed algebra together
with other $5$ Noether-like operators. \newline
\newline
\textbf{Case }$\mathcal{N}_{3,6}^{1}$\newline
We obtain six Noether-like operators by using (\ref{eq6})
\begin{eqnarray}
\mathbf{X}_{1} =\frac{\partial }{\partial f},\quad \mathbf{X}_{2}=\frac{\partial
}{\partial g},\quad \mathbf{X}_{3}=x\frac{\partial }{\partial x}+f\frac{\partial }{%
\partial f}+g\frac{\partial }{\partial g}, \nonumber \\
 \mathbf{X}_{4}=g\frac{\partial }{%
\partial f}-f\frac{\partial }{\partial g},  \quad 
\mathbf{X}_{5} =2xf\frac{\partial }{\partial x}+(f^{2}-g^{2})\frac{%
\partial }{\partial f}+2fg\frac{\partial }{\partial g},  \nonumber \\
\mathbf{X}_{6} =2xg\frac{\partial }{\partial x}+2fg\frac{\partial }{%
\partial f}-(f^{2}-g^{2})\frac{\partial }{\partial g}, \label{eq71}
\end{eqnarray}%
for the following system of ODEs%
\begin{eqnarray}
xf^{\prime \prime } =-f^{\prime 3}+3f^{\prime }g^{\prime 2}-f^{\prime }/2,
\nonumber \\
xg^{\prime \prime } =-3f^{\prime 2}g^{\prime }+g^{\prime 3}-g^{\prime }/2, \label{eq72}
\end{eqnarray}%
which arise from the r-CODE in Table 1. By the utilization of (\ref{eq4}) it can be verified that the above system admits
the Lagrangians
\begin{eqnarray}
L_{1} =\frac{f^{\prime }}{x}-\frac{f^{\prime }}{2x(f^{\prime 2}+g^{\prime
2})},  \quad
L_{2} =\frac{g^{\prime }}{x}+\frac{g^{\prime }}{2x(f^{\prime 2}+g^{\prime
2})}. \label{eq73}
\end{eqnarray}%
Our next step is to construct first integrals for the system (\ref{eq72})
by taking into account all Noether-like operators. Now starting with
the Noether-like operators $\mathbf{X}_{1}$ and $\mathbf{X}_{2}$ we
arrive at the two real first integrals%
\begin{eqnarray}
I_{1} =\frac{1}{x}+\frac{f^{^{\prime }2}-g^{\prime 2}}{%
2x(f^{^{\prime }2}+g^{\prime 2})^{2}} ,  \quad
I_{2} =\frac{f^{\prime }g^{\prime
}}{x(f^{^{\prime }2}+g^{\prime 2})^{2}} .\label{eq74}
\end{eqnarray}%
The other first integrals corresponding to $\mathbf{X}_{3}$ and $\mathbf{X}%
_{4}$ are
\begin{eqnarray}
I_{3} =-\frac{f^{\prime
}}{f^{2}+g^{2}}+\frac{f}{x}+ \frac{f(f^{^{\prime
}2}-g^{\prime 2})+2gf^{\prime }g^{\prime }}{2x(f^{^{\prime
}2}+g^{\prime 2})^{2}} ,  \nonumber \\
I_{4} =\frac{g^{\prime
}}{f^{2}+g^{2}}+\frac{g}{x}+ \frac{g(f^{^{\prime
}2}-g^{\prime 2})-2ff^{\prime }g^{\prime }}{2x(f^{^{\prime
}2}+g^{\prime 2})^{2}} . \label{eq75}
\end{eqnarray}
Similarly, $\mathbf{X}_{5}$ and $\mathbf{X}_{6}$ yield
\begin{eqnarray}
I_{5} =2x-\frac{(f^{2}-g^{2})}{x}+\frac{2(ff^{\prime }+gg^{\prime })}{%
f^{^{\prime }2}+g^{^{\prime }2}}+\frac{(f^{2}-g^{2})(f^{\prime 2}-g^{\prime
2})+4ff^{\prime }gg^{\prime }}{2x(f^{^{\prime }2}+g^{^{\prime }2})^{2}},
\nonumber \\
I_{6} =\frac{2fg}{x}-\frac{2(f^{\prime }g-fg^{\prime })}{f^{^{\prime
}2}+g^{^{\prime }2}}+\frac{fg(f^{\prime 2}-g^{\prime 2})-f^{\prime
}g^{\prime }(f^{2}-g^{2})}{x(f^{^{\prime }2}+g^{^{\prime }2})^{2}}. \label{eq76}
\end{eqnarray}%
The closed algebra of these Noether-like operators is spanned by%
\begin{eqnarray}
\lbrack \mathbf{X}_{1},\mathbf{X}_{2}] =[\mathbf{X}_{3},\mathbf{X}_{4}]=[%
\mathbf{X}_{5},\mathbf{X}_{6}]=0,  \nonumber \\
\lbrack \mathbf{X}_{1},\mathbf{X}_{3}] =\mathbf{X}_{1},\quad [\mathbf{X}%
_{1},\mathbf{X}_{4}]=-\mathbf{X}_{2},\quad [\mathbf{X}_{1},\mathbf{X}%
_{5}]=2\mathbf{X}_{3},  \nonumber \\
\lbrack \mathbf{X}_{1},\mathbf{X}_{6}] =2\mathbf{X}_{4},\quad [\mathbf{X}%
_{2},\mathbf{X}_{3}]=\mathbf{X}_{2},\quad[\mathbf{X}_{2},\mathbf{X}_{4}]=%
\mathbf{X}_{1},  \nonumber \\
\lbrack \mathbf{X}_{2},\mathbf{X}_{5}] =-2\mathbf{X}_{4},\quad[\mathbf{X%
}_{2},\mathbf{X}_{6}]=2\mathbf{X}_{3},\quad [\mathbf{X}_{3},\mathbf{X}%
_{5}]=\mathbf{X}_{5},  \nonumber \\
\lbrack \mathbf{X}_{3},\mathbf{X}_{6}] =\mathbf{X}_{6},\quad [\mathbf{X}%
_{4},\mathbf{X}_{5}]=\mathbf{X}_{6},\quad [\mathbf{X}_{4},\mathbf{X}_{6}]=-\mathbf{%
X}_{5}. \label{eq77}
\end{eqnarray}%
Indeed these operators span a six-dimensional algebra.
\newline
\newline
\textbf{Remaining Cases:}\newline
\textbf{Case }$\mathcal{N}_{3,5}^{3}$\newline
Consider the pair of two real Lagrangians
\begin{eqnarray}
L_{1} =\frac{-r^{1/2}}{(1+b^{2})}e^{\frac{b}{2}\theta _{1}}\left (\cos
(\theta/2)\right )+xf^{\prime },  \nonumber \\
L_{2} =\frac{-r^{1/2}}{(1+b^{2})}e^{\frac{b}{2}\theta _{1}}\left (\sin
(\theta /2)\right)+xg^{\prime }. \label{eq78}
\end{eqnarray}%
The system of second-order ODEs associated with the above
Lagrangians is
\begin{eqnarray}
f^{\prime \prime } =r^{3/2}e^{-\frac{b}{2}\theta _{1}}\left(\cos
(3\theta /2)\right),\quad
g^{\prime \prime } =r^{3/2}e^{-\frac{b}{2}\theta _{1}}\left(\sin
(3\theta /2)\right), \label{eq79}
\end{eqnarray}%
where $r$ and $\theta $ are given by (\ref{eq59}) and $\theta _{1}$ by
(\ref{eq64}). These correspond to the r-CODE in Table 1. The three
Noether-like operators admitted by (\ref{eq79}) are
\begin{equation}
\mathbf{X}_{1}=\frac{\partial }{\partial x},\quad \mathbf{X}_{2}=\frac{\partial }{%
\partial f},\quad \mathbf{X}_{3}=\frac{\partial }{\partial g}. \label{eq80}
\end{equation}%
The first integrals are
\begin{eqnarray}
I_{1} =x+\frac{r^{-1/2}}{(1+b^{2})}e^{\frac{b}{2}\theta
_{1}}\left((b-f^{\prime
})\cos (\theta /2)-g^{\prime }\sin (\theta /2)\right),  \nonumber \\
I_{2} =-\frac{r^{-1/2}}{(1+b^{2})}e^{\frac{b}{2}\theta
_{1}}\left ((b-f^{\prime })\sin (\theta /2)+g^{\prime }\cos (\theta /2)\right), \label{eq81}
\end{eqnarray}%
for the system (\ref{eq79}) associated with $\mathbf{X}_{2}$ and
$\mathbf{X}_{3}.$ Now corresponding to $\mathbf{X}_{1}$ we have
\begin{eqnarray}
I_{3} =\frac{r^{-1/2}}{(1+b^{2})}e^{\frac{b}{2}\theta
_{1}}\left((1+bf^{\prime
})\cos (\theta /2)+bg^{\prime }\sin (\theta /2)\right)+f,  \nonumber \\
I_{4} =\frac{r^{-1/2}}{(1+b^{2})}e^{\frac{b}{2}\theta
_{1}}\left(bg^{\prime }\cos (\theta /2)-(1+bf^{\prime })\sin (\theta
/2)\right)+g. \label{eq82}
\end{eqnarray}
We can obtain two more first integrals corresponding to other Noether operators. These
operators correspond to alternate Lagrangians. These integrals can be expressed as a combination of
above four first integrals. \newline
\newline
\textbf{Case }$\mathcal{N}_{3,6}^{2}$\newline
Consider the system of two second-order ODEs%
\begin{eqnarray}
xf^{\prime \prime } =f^{\prime ^{3}}-3f^{\prime }g^{\prime 2}+f^{\prime
}+r_{1}^{3/2}\left (A\cos (3\theta_{1} /2)\right),  \nonumber \\
xg^{\prime \prime } =3f^{\prime ^{2}}g^{\prime }-g^{\prime 3}+g^{\prime
}+r_{1}^{3/2}\left(A\sin (3\theta_{1} /2)\right), \label{eq83}
\end{eqnarray}%
where $r_{1}$ and $\theta_{1} $ are given by (\ref{eq64}). This system
arises from the r-CODE in Table 1.
 The Lagrangians of the above
system are readily found to be of the forms
\begin{eqnarray}
L_{1} =\frac{r_{1}^{1/2}}{x}\cos (\theta_{1} /2)+\frac{Af^{\prime }}{x},  \quad
L_{2} =\frac{r_{1}^{1/2}}{x}\sin (\theta_{1} /2)+\frac{Ag^{\prime }}{x}. \label{eq84}
\end{eqnarray}%
The six Noether-like operators related to (\ref{eq83}) are shown to be%
\begin{eqnarray}
\mathbf{X}_{1} =\frac{\partial }{\partial f},\quad \mathbf{X}_{2}=\frac{\partial
}{\partial g},\quad 
\mathbf{X}_{3} =x\frac{\partial }{\partial x}+f\frac{\partial }{\partial f}%
+g\frac{\partial }{\partial g},\nonumber \\
 \mathbf{X}_{4}=g\frac{\partial }{\partial f}-f%
\frac{\partial }{\partial g},  \quad
\mathbf{X}_{5} =2xf\frac{\partial }{\partial x}+(f^{2}-g^{2}-x^{2})\frac{%
\partial }{\partial f}+2fg\frac{\partial }{\partial g},  \nonumber \\
\mathbf{X}_{6} =2xg\frac{\partial }{\partial x}+2fg\frac{\partial }{%
\partial f}-(f^{2}-g^{2}-x^{2})\frac{\partial }{\partial g}. \label{eq85}
\end{eqnarray}%
We take $\mathbf{X}_{1}$ and $\mathbf{X}_{2}$ which imply the first
integrals
\begin{eqnarray}
I_{1} =\frac{A}{x}+\frac{r_{1}^{-1/2}}{x}\left(f^{\prime }\cos (\theta_{1}
/2)+g^{\prime }\sin (\theta_{1} /2)\right),  \nonumber \\
I_{2} =\frac{1}{x}r_{1}^{-1/2}\left(f^{\prime }\sin (\theta_{1} /2)-g^{\prime }\cos
(\theta_{1} /2)\right), \label{eq86}
\end{eqnarray}%
for (\ref{eq83}). Other Noether-like operators $\mathbf{X}_{3}$ and
$\mathbf{X}_{4}$ enables us to write the first integrals
\begin{eqnarray}
I_{3} =r_{1}^{-1/2}\left(\left(1+\frac{ff^{\prime }-gg^{\prime }}{x}\right)\cos\frac{\theta_{1}}{2}+%
\frac{(f^{\prime }g+fg^{\prime })}{x}\sin \frac{\theta_{1}}{2}\right)+\frac{Af}{x},  \nonumber
\\
I_{4} =r_{1}^{-1/2}\left(\frac{(f^{\prime }g+fg^{\prime })}{x}\cos \frac{\theta_{1}}{2}-\left(1+%
\frac{ff^{\prime }-gg^{\prime }}{x}\right)\sin \frac{\theta_{1}}{2}\right)+\frac{Ag}{x}, \label{eq87}
\end{eqnarray}%
The remaining two invariants related to the Noether-like operators $\mathbf{X%
}_{5}$, $\mathbf{X}_{6}$ are found to be
\begin{eqnarray}
I_{5} =A\left (x+\frac{f^{2}-g^{2}}{x}\right )-r_{1}^{-1/2}%
\left (xf^{\prime }-\frac{(f^{2}-g^{2})f^{\prime }-2fgg^{\prime }}{x}\right )\cos \frac{\theta_{1} }{2}
\nonumber \\
-r_{1}^{-1/2}\left (xg^{\prime }-2f-2g+\frac{(f^{2}-g^{2})g^{\prime
}+2fgf^{\prime }}{x}\right )\sin \frac{\theta_{1} }{2},  \nonumber \\
I_{6} =\frac{2Afg}{x}+r_{1}^{-1/2}\left (xf^{\prime }-\frac{%
f^{\prime }(f^{2}-g^{2})-2fgg^{\prime }}{x}\right)\sin \frac{\theta_{1} }{2}  \nonumber \\
+r_{1}^{-1/2}\left(xg^{\prime }-2f-2g-\frac{g^{\prime
}(f^{2}-g^{2})+2f^{\prime }fg}{x}\right)\cos \frac{\theta_{1} }{2}. \label{eq88}
\end{eqnarray}
\newline
\textbf{Case} $\mathcal{N}_{3,6}^{3}$\newline
Lastly, a nonlinear system of ODEs can be mapped to the the r-CODE%
\begin{equation}
u^{\prime \prime }=A\left( \frac{1+u^{\prime 2}+(u-xu^{\prime })^{2}}{%
1+x^{2}+u^{2}}\right) ^{3/2}, \label{eq89}
\end{equation}%
which admits the $3-$dimensional complex Noether algebra $\mathbf{Z}%
_{1}=(1+x^{2})\partial_{x}+xu\partial_{u}$, $\mathbf{Z}%
_{2}=xu\partial_{x}+(1+u^{2})\partial_{u}$, $\mathbf{Z}%
_{3}=u\partial_{x}-x\partial_{u}$. By choosing $\mathbf{Z}%
_{1},$ we have the complex first integral%
\begin{equation}
I_{1}=A\frac{\alpha }{\sqrt{1+\alpha ^{2}}}+\frac{1}{\sqrt{1+\alpha
^{2}+\beta ^{2}}}, \label{eq90}
\end{equation}%
where $\alpha =u(1+x^{2})^{-1/2}$, $\beta =u^{\prime
}(1+x^{2})^{1/2}-xu(1+x^{2})^{-1/2}.$ The complex Lagrangian of the above
r-CODE is
\begin{equation}
L=\frac{1}{t^{2}}(1+x^{2})^{-3/2}\left( (\gamma -u^{\prime })\sin \theta
+\delta \sec \theta \right) +d(x,u)u^{\prime }+e(x,u), \label{eq91}
\end{equation}%
where $d$ and $e$ satisfy $d_{x}=e_{u}+A(1+x^{2}+u^{2})^{-3/2}$ and $t=\sqrt{%
1+x^{2}+u^{2}}/(1+x^{2})$, $\delta =tu,\tan \theta =(u^{\prime }-\gamma
)/t$ and $\gamma =xu/(x^{2}+1).$ The other two complex first
integrals
corresponding to the respective $\mathbf{Z}_{2}$ and $\mathbf{Z}_{3}$ are given as%
\begin{eqnarray}
I_{2} =\frac{Ax}{\sqrt{1+x^{2}+u^{2}}}-\frac{u^{\prime }}{\sqrt{%
1+u^{\prime 2}+(u-xu^{\prime })^{2}}},  \nonumber \\
I_{3} =\frac{A}{\sqrt{1+x^{2}+u^{2}}}+\frac{(x^{2}+1)^{-3/2}}{t^{2}}\left(
\delta\sin \theta +\delta \sec \theta \right) , \label{eq92}
\end{eqnarray}%
where $t$, $\gamma$, $\delta$ and $\theta$ are given above and $\delta
=x+u\gamma .$ The complex Lagrangian, Noether symmetries and first
integrals exactly correspond to the emerging system of (\ref{eq89}).
\section{Maximal Noether algebra for two-dimensional systems}

For many problems of physical interest, linear systems as well as
some nonlinear systems, the Lagrangians that arise do result in the
maximum number of Noether point symmetries for the said system, are
known as standard Lagrangians. These Lagrangians are constructed in
a straightforward and obvious manner from physical considerations. A
Lagrangian of the given dynamical equation is said to be maximally
symmetric if it gives the maximum number of Noether symmetries that
can occur for the dynamical equation. We know that for a dynamical
equation, infinitely many Lagrangians may exist, but regarding one
Lagrangian, there is a unique dynamical equation. The symmetries of
Lagrangian admitted by the systems of ODEs is of great significance
especially in regard to their use in physical applications. Each
Noether symmetry for the given system provides one constant of
motion. The analytic continuation furnishes two nice invariant
quantities for the corresponding system of ODEs. The conservation
laws play a vital role in the study of physical phenomena. As an
example it would be significant to mention that if a dynamical
system remains invariant under translation in time then it has a
conservation of energy. Likewise, a rotational symmetry correspond
to the conservation of angular momentum. The complex symmetries
correspond to conformal motions. Therefore, if a differential equation and its
Lagrangian admit a complex rescaling symmetry then that dynamical
system respects conservation of infinitesimal angles.

We now focus our attention on a two-dimensional system of free particle
equations. The determination of Noether symmetries and all invariants for
this simplest system
\begin{eqnarray}
f^{\prime \prime } =0,\quad 
g^{\prime \prime } =0,\label{eq93}
\end{eqnarray}%
can be obtained via the classical Noether approach. We utilize the complex approach here.
 The maximal nine
Noether-like operators for (\ref{eq93}) relative to the Lagrangians $L_{1}=\frac{1%
}{2}(f^{\prime 2}-g^{\prime 2})$ and $L_{2}=f^{\prime }g^{\prime }$ are
found to be (note that these are not standard)%
\begin{eqnarray}
\mathbf{X}_{1} =\frac{\partial }{\partial x},\quad \mathbf{X}_{2}=\frac{\partial
}{\partial f},\quad \mathbf{X}_{3}=\frac{\partial }{\partial g},\quad \mathbf{X}_{4}=2x%
\frac{\partial }{\partial x}+f\frac{\partial }{\partial f}+g\frac{\partial }{%
\partial g},  \nonumber \\
\mathbf{X}_{5} =g\frac{\partial }{\partial f}-f\frac{\partial }{\partial g}%
,\quad \mathbf{X}_{6}=x\frac{\partial }{\partial f},\quad \mathbf{X}_{7}=x\frac{\partial
}{\partial g}, \nonumber \\
\mathbf{X}_{8}=x^{2}\frac{\partial }{\partial x}+x(f\frac{%
\partial }{\partial f}+g\frac{\partial }{\partial g}),\quad 
\mathbf{X}_{9} =x(g\frac{\partial }{\partial f}-f\frac{\partial }{\partial
g}). \label{eq94}
\end{eqnarray}%
Notice that the operators $\mathbf{X}_{1},...,\mathbf{X}_{8}$ are
the usual Noether symmetries with respect to the standard Lagrangian
of (\ref{eq93}). Surprisingly, the maximal $8-$dimensional Noether algebra
(excluding $\mathbf{X}_{9}$) for the above system can be obtained
from the analytic continuation of the maximal $5-$dimensional
Noether algebra, i.e.,
\begin{eqnarray}
\mathbf{Z}_{1}=\frac{\partial }{\partial x},\quad \mathbf{Z}_{2}=\frac{\partial }{%
\partial u},\quad \mathbf{Z}_{3}=2x\frac{\partial }{\partial x}+u\frac{\partial }{%
\partial u},\quad \mathbf{Z}_{5}=x\frac{\partial }{\partial u},\nonumber \\
\mathbf{Z}_{4}=x^{2}\frac{\partial }{\partial x}+xu\frac{%
\partial }{\partial u}, \label{eq95}
\end{eqnarray}%
for the equation $u^{\prime \prime }=0,$ relative to the Lagrangian 
$L=u^{\prime 2}/2,$ (termed as the standard Lagrangian since it is the kinetic
energy of the particle) which is maximally symmetric since it yields the
maximum number of Noether point symmetries (\ref{eq95}). The extra operator $\mathbf{X}%
_{9}$ seems to be a symmetry corresponding to $\mathbf{Z}_{4}$. It
does not satisfy the classical Lie conditions and it thus fails to
be a symmetry of (\ref{eq93}). Yet it \textit{yields} first integrals, see
$I_{4,1}$ and $I_{4,2},$ in Table 2. Hence, this surplus symmetry is
of an extraordinary characteristic. It is, somehow, attached to the
system in a very strange way, although it is not a Lie symmetry. Indeed this 
extra ordinary feature of complex variables plays a central role in the 
determination of invariants of two-dimensional systems thereby yields 
deeper insights into the analysis of physical models (\cite{28}). \newline\newline
Table 2: First Integrals for simplest two dimensional system\newline
\begin{tabular}{|l|l|l|}
\hline
Symmetries  ~~~~Complex First Integrals  ~~~~First Integrals \\ \hline
$\partial _{x}$ ~~~~~~~~~~~~~~~~ $I_{1}=u^{\prime 2}$  $%
\begin{tabular}{l}
~~~~~~~~~~~~~~~~~~~~~~$I_{1,1}=f^{\prime 2}-g^{\prime 2}$ \\
~~~~~~~~~~~~~~~~~~~~~~$I_{1,2}=f^{\prime }g^{\prime }$%
\end{tabular}%
$ \\ \hline
$\partial _{u}$  ~~~~~~~~~~~~~~~~~$I_{2}=u^{\prime }$
\begin{tabular}{l}
~~~~~~~~~~~~~~~~~~~~~~~$I_{2,1}=f^{\prime }$ \\
~~~~~~~~~~~~~~~~~~~~~~~$I_{2,2}=g^{\prime }$%
\end{tabular}
\\ \hline
$x\partial _{u}$ ~~~~~~~~~~~~~~~$I_{3}=xu^{\prime }-u$
\begin{tabular}{l}
~~~~~~~~~~~~~~~~~$I_{3,1}=xf^{\prime }-f$ \\
~~~~~~~~~~~~~~~~~$I_{3,2}=xg^{\prime }-g$%
\end{tabular}
\\ \hline
$2x\partial _{x}+u\partial _{u}$  ~~~~~~$I_{4}=-xu^{\prime 2}+uu^{\prime }$
\begin{tabular}{l}
~~~~~~~~~~$I_{4,1}=-x(f^{\prime 2}-g^{\prime 2})+ff^{\prime }-gg^{\prime }$ \\
~~~~~~~~~~$I_{4,2}=-xf^{\prime }g^{\prime }+fg^{\prime }+f^{\prime }g$%
\end{tabular}
\\ \hline
$x^{2}\partial _{x}+xu\partial _{u}$  ~~~~$I_{5}=xuu^{\prime }-\frac{%
x^{2}u^{\prime 2}}{2}-\frac{u^{2}}{2}$
\begin{tabular}{l}
~~~~$I_{5,1}=x(ff^{\prime }-gg^{\prime })-\frac{x^{2}}{2}(f^{\prime 2}-g^{\prime
2})-$ \\
~~~~$\frac{1}{2}(f^{2}-g^{2})$ \\
~~~~$I_{5,2}=x(fg^{\prime }+f^{\prime }g)-x^2f^{\prime }g^{\prime }-fg$%
\end{tabular}
\\ \hline
\end{tabular}%
\newline
\newline

We know that for any particle Lagrangian, the maximum dimension of the
Noether algebra is five \cite{18}. As mentioned earlier, the free particle equation
admits another Lagrangian of the type $L=u^{2}u^{\prime }\log u^{\prime
}-u^{2}u^{\prime }-xuu^{\prime 2}+\frac{1}{6}x^{2}u^{\prime 3}+au^{\prime
}+b,$ where $a$ and $b$ are related by $a_{x}=b_{u}.$ The Noether symmetries
of this Lagrangian are
\begin{equation}
\mathbf{Z}_{6}=u\frac{\partial }{\partial x},\quad \mathbf{Z}_{7}=x\frac{\partial
}{\partial x},\quad \mathbf{Z}_{8}=xu\frac{\partial }{\partial x}+u^{2}\frac{%
\partial }{\partial u}. \label{eq96}
\end{equation}%
The Noether symmetries of the two Lagrangians together give the
complete symmetry generators of $sl(3,\mathbb{C})$ for the restricted
particle equation. In the solution of simple harmonic oscillator the
two Lagrangian pictures are essential in order to describe fully the
algebra and periodicity of motion of the oscillator \cite{27}. The
complete algebraic description, i.e., the complete determination of
all the physical constants is obtained from these. We also know that
the simplest system admits $sl(4,\Re).$ The question of complete
determination of invariants of scalar linear equations can be
extended to system of two free particle equations and it can be
asked whether one can find the rest of the first integrals
corresponding to a $7-$dimensional Noether subalgebra of $sl(4,\Re).$
Therefore, we present  Conjecture 1.
\newline
\textbf{Conjecture 1.} \textit{The $7-$dimensional subalgebra of $sl(4,\Re)$ for the simplest system which complements the
$8-$dimensional Noether algebra is derivable from a complex subalgebra of $5-$%
dimension of an r-CODE via an alternative Lagrangian}. \newline
It is important to mention that the use of alternative Lagrangians grants
more Noether-like operators which would determine complete set of physical
constants for the simplest system. We know that every linear second-order
r-CODE%
\begin{equation}
u^{\prime \prime }=\alpha (x)u^{\prime }+\beta (x)u+\gamma (x), \label{eq97}
\end{equation}%
is equivalent to the simplest r-CODE, $U^{\prime \prime }=0$ via an
appropriate point transformation \cite{1}. It is also known that it
admits a $5-$dimensional Noether algebra with respect to a standard
Lagrangian. Thus, for systems of linear ODEs which can
be mapped to a linear r-CODE the Noether-like operators may induce an $8-$%
dimensional Noether algebra. Consequently, we state the Conjecture $2$.
\newline
\textbf{Conjecture 2.} \textit{Every system of two linear second-order ODEs
that is obtainable from a linear r-CODE with usual Lagrangian admits the 
$8-$dimensional Noether algebra}.

\section{Conclusion}

In this paper, we have addressed the problem of Noether classification
 for systems of two second-order EL ODEs that arise from the submaximal and
 maximal Noether symmetry classification of Lagrangians on the line. We achieved this by
 the introduction of Noether-like operators. In this way, we provided an algebraic
study of systems of ODEs which are variational and obtainable from
complex EL equations. Also this study provides new ways of obtaining
first integrals for those classes of systems of nonlinear equations
that have not been studied before via the classical Noether
approach. It is certainly meaningful to classify all such systems
which exhibit the structure of Noether-like operators. We have
discussed these for five, six and eight Noether-like operators. In
some cases the Noether-like operators turn out to be the Noether
point symmetries. The Noether counting theorem, an analogue of the
Lie
counting theorem, states that a second-order  EL r-CODE can have $0,1,2,3$ or $5-$%
dimensional Noether algebra \cite{18}. We have shown that a system
of two second-order ODEs that
arise from an r-CODE can have $%
0,1,2,3,4,5,6$ or $9-$dimensional Noether-like operators. The lower dimensional cases up to dimension 4 were illustrated by means
of examples. The simplest system
possesses the maximal $8-$dimensional Noether point symmetry algebra which we achieved
remarkably by the analytic continuation of the $5-$dimensional Noether algebra
of a linear second-order r-CODE. In addition it supplies an extra operator which is not a
Noether symmetry of the system, although it does permit the emergence of a
first integral.\newline

{\bf Acknowledgements}. FM thanks the Higher Education Commission (HEC) of Pakistan for a visiting
professorship during which this research was commenced.

\end{document}